\documentclass{amsart}
\usepackage{amssymb,latexsym,amsmath,amsfonts,hhline}
\usepackage{pdfsync}
\usepackage{xcolor} 
\usepackage{comment,colortbl}
\input xy
\xyoption{all}

\newcommand{\cal}{\mathcal}

\newtheorem{formula}{}[section]
\newtheorem{definition}[formula]{Definition}
\newtheorem{corollary}[formula]{Corollary}
\newtheorem{remark}[formula]{Remark}
\newtheorem{lemma}[formula]{Lemma}
\newtheorem{theorem}[formula]{Theorem}
\newtheorem{example}[formula]{{\bf Example}}

\def\thrm{\begin{theorem}}
\def\thrml#1{\begin{theorem}\label{#1}}
\def\ethrm{\end{theorem}}
\def\rmrk{\begin{remark}}
\def\rmrkl#1{\begin{remark}\label{#1}}
\def\ermrk{\end{remark}}
\def\dfntn{\begin{definition}}
\def\dfntnl#1{\begin{definition}\label{#1}}
\def\edfntn{\end{definition}}
\def\nmrt{\begin{enumerate}}
\def\enmrt{\end{enumerate}}
\def\tm#1{\item[{\rm (#1)}]}
\def\qtn{\begin{equation}}
\def\qtnl#1{\begin{equation}\label{#1}}
\def\eqtn{\end{equation}}
\def\lmm{\begin{lemma}}
\def\lmml#1{\begin{lemma}\label{#1}}
\def\elmm{\end{lemma}}
\def\crllr{\begin{corollary}}
\def\crllrl#1{\begin{corollary}\label{#1}}
\def\ecrllr{\end{corollary}}
\def\css{\begin{cases}}
\def\ecss{\end{cases}}
\def\xmpl{\begin{example}}
\def\xmpll#1{\begin{example}\label{#1}}
\def\exmpl{\end{example}}

\def\proof{\noindent{\bf Proof}.\ }

\def\cA{{\cal A}}

\def\cP{{\cal P}}

\def\cX{{\cal X}}

\def\fX{{\frak X}}

\DeclareMathOperator{\aut}{Aut}

\DeclareMathOperator{\dom}{Dom}

\DeclareMathOperator{\id}{id}

\DeclareMathOperator{\iso}{Iso}

\DeclareMathOperator{\sym}{Sym}

\def\ddp{{\scriptscriptstyle\bot}}

\def\qaq{\quad\text{and}\quad}

\def\wh{\widehat}
\def\wt{\widetilde}

\newcommand{\grp}[1]{\langle {#1}\rangle}

\def\VRT#1{*=<5mm>[o][F-]{#1}}

\def\grphp#1{$\xymatrix@R=10pt@C=10pt@M=0pt@L=2pt{#1}$}

\def\proof{{\bf Proof}.\ }
\def\eprf{\hfill$\square$}

\begin{document}

\title{Two-valenced association schemes and the Desargues theorem}
\author{Mitsugu Hirasaka}
\address{Department of Mathematics, Pusan National University, Jang-jeon dong, Busan, Republic of Korea}
\email{hirasaka@pusan.ac.kr}
\thanks{This research was supported by Basic Science Research Program through the National Research Foundation of Korea(NRF) funded by the Ministry of  Science, ICT \& Future Planning(NRF-2016R1A2B4013474).}
\author{Kijung Kim}
\address{Department of Mathematics, Pusan National University, Jang-jeon dong, Busan, Republic of Korea}
\email{knukkj@pusan.ac.kr}
\thanks{The second author was supported by Basic Science Research Program through the National Research Foundation of Korea funded by the Ministry of Education (NRF-2017R1D1A3B03031349)}
\author{Ilia Ponomarenko}
\address{St.Petersburg Department of the Steklov Mathematical Institute, St.Petersburg, Russia}
\email{inp@pdmi.ras.ru}
\thanks{The work of the third author was supported by the RAS Program of Fundamental Research ``Modern Problems of Theoretical Mathematics''. }

\begin{abstract}
The main goal of the paper is to establish a sufficient condition for a two-valenced association scheme to be schurian and separable. To this end, an analog of  the Desargues theorem is introduced for a  noncommutative geometry defined by the scheme in question. It turns out that if the geometry has enough many Desarguesian configurations, then under a technical condition the scheme is schurian and separable. This result enables us to give short proofs for known statements on the schurity and separability of quasi-thin and pseudocyclic schemes. Moreover, by the same technique we prove a new result: given a prime~$p$, any  $\{1,p\}$-scheme with thin residue isomorphic to an elementary abelian $p$-group of rank greater than two, is schurian and separable.
\end{abstract}

\maketitle

\section{Introduction}

One of the fundamental problems in theory of association schemes is to determine whether a given scheme is {\it schurian}, i.e., comes from a permutation group, and/or {\it separable}, i.e., uniquely determined by its intersection number array (for the exact definitions, see Section~\ref{200817a}). In the two last decades these two problems are intensively studied for the {\it two-valenced} schemes see, e.g., ~\cite{Hirasaka2001,Muzychuk2008a,Muzychuk2012a,ChenPonom2017}; here an association scheme is said to be {\it two-valenced} if the valencies of its basic relations take exactly two values, and if  they are $1$ and $k$, the term {\it $\{1,k\}$-valenced} scheme is also used.\medskip

An analysis of the known proofs that certain two-valenced schemes are schurian or separable shows that in all cases the following two properties are significant. The first one is that there are sufficiently many intersection numbers of the scheme in question that are equal to~$1$; to define this property precisely, we introduce in Section~\ref{281117l} the {\it saturation} condition (a special case of it appeared in~\cite{ChenPonom2017}). The second property expresses the fact that in a noncommutative ``affine'' geometry determined by the two-valenced scheme, there are sufficiently many  Desarguesian configurations, see Section~\ref{281117a4}.\medskip

A two-valenced scheme having the first and second properties is said to be  {\it saturated} and {\it Desarguesian}, respectively. A model example illustrating these two properties is the scheme of a finite affine space that is two-valenced, saturated (except for few cases), and Desarguesian if the dimension of the affine space is at least~$3$ (see Examples in Sections~\ref{281117l} and~\ref{281117a4}).


\thrml{181117a}
Let $\cX$ be a two-valenced scheme. Assume that $\cX$ is saturated and  Desarguesian. Then $\cX$ is schurian and  separable.
\ethrm

For the scheme of a finite affine space  of dimension  at least~$3$, Theorem~\ref{181117a} expresses a well-known fact that in the Desarguesian case, this scheme is reconstructed  from its automorphism group and is uniquely determined by its order.\medskip

The power of Theorem~\ref{181117a} is illustrated by the statements below. The first two of them are known results, but using the theorem we are able to give much shorter and clear proofs than in the original papers. 

\crllrl{110317a}{\rm \cite{Muzychuk2012a}}
There exists a function $f$ such that any pseudocyclic scheme of valency $k>1$ and degree at least $f(k)$ is schurian and separable.
\ecrllr

From the proof of Corollary~\ref {110317a} given in Section~\ref{090818a}, it follows that the function~$f$ satisfies the inequality $f(k)\le 1+ 3k^6$. A more subtle arguments used in~\cite{ChenPonom2017} show that $f(k)\le 1+6k(k-1)^2$. 

\crllrl{280618a}{\rm \cite{Muzychuk2008a,Muzychuk2012}}
Any quasi-thin scheme satisfying the condition $n_{s^*s}\ne 2$ for all basis relations~$s$, is schurian and separable.
\ecrllr

Except for the separability statement, the result of Corollary~\ref{280618a} is contained in~\cite{Muzychuk2008a}. On the other hand, the schurian and separable quasi-thin schemes were characterized in~\cite{Muzychuk2012}. This result can also be deduced from an analog of Theorem~\ref{181117a}, in which the Desarguesian condition is replaced by a weaker one (namely, the amount of the required Desarguesian configuration is reduced). To keep the text more compact, we do not go into detailed explanation of this topic.\medskip

The second main result of the present paper concerns the schurity and separability of a class of meta-thin schemes introduced in~\cite{Hirasaka2005}. A meta-thin scheme can be thought as an extension of a regular scheme $\cX$ by another regular scheme. Even if the meta-thin scheme is a $\{1,p\}$-scheme, where $p$ is a prime, the schurity and separability problems seem to  be very complicated. In a sense, the answer depends on the scheme~$\cX$, which can be chosen as the thin residue of the scheme in question. For example, if the group associated with $\cX$ has distributive lattice  of normal subgroups, then the scheme in question is schurian and separable~\cite{Hirasaka2018}. However, there are many non-schurian and non-separable meta-thin $\{1,p\}$-schemes for which that group is elementary abelian of order~$p^2$~\cite{Hirasaka2016}. The following theorem, which is also deduced from Theorem~\ref{181117a}, shows that if the group is an elementary abelian $p$-group of rank greater than two, then the situation is smooth.

\thrml{thm111}
Given a prime $p$, any  $\{1,p\}$-scheme with thin residue isomorphic to an elementary abelian $p$-group of rank greater than two, is schurian and separable.
\ethrm

The paper is organized as follows. For the reader convenience a background on association schemes  and related concepts is given in Section~\ref{200817a}. In Sections~\ref{281117l} and~\ref{281117a4} we introduce and study the saturated and Desarguesian two-valenced schemes, respectively. The proof of Theorem~\ref {181117a} is given in Section~\ref{090818n}. In the final Section~\ref{090818a} we prove the corollaries and Theorem~\ref{thm111}.\pagebreak

{\bf Notation.}

Throughout the paper, $\Omega$ denotes a finite set.

The diagonal of the Cartesian product $\Omega\times\Omega$ is denoted by~$1_\Omega$ or~$1$.

For a relation $r\subseteq\Omega\times\Omega$, we set $r^*=\{(\beta,\alpha):\ (\alpha,\beta)\in r\}$ and $\alpha r=\{\beta\in\Omega:\ (\alpha,\beta)\in r\}$ for all $\alpha\in\Omega$.

For a relation $r\subseteq\Omega\times\Omega$ and sets $\Delta,\Gamma\subseteq\Omega$, we set $s_{\Delta,\Gamma}=s\cap(\Delta\times\Gamma)$. If $S$ is a set of relations, we put $S_{\Delta,\Gamma}=\{s_{\Delta,\Gamma}:\ s\in S\}$ and denote $S_{\Delta,\Delta}$ by $S_\Delta$.

For relations $r,s\subseteq \Omega\times\Omega$, we set
$r\cdot s=\{(\alpha,\beta):\ (\alpha,\gamma)\in r,\ (\gamma,\beta)\in s$
for some $\gamma\in\Omega\}$. If $S$ and $T$ are sets of relations, we set
$S\cdot T=\{s\cdot t:\ s\in S,\, t\in T\}$.

For a set $S$ of relations on~$\Omega$, we denote by $S^\cup$ the set of all unions of the elements of $S$, and put $S^*=\{s^*:\ s\in S\}$
and $\alpha S=\cup_{s\in S}\alpha s$, where $\alpha\in\Omega$.

An elementary abelian $p$-group of order $p^m$ is denoted by $E_{p^m}$.

\section{Association schemes}\label{200817a}

In our presentation of association schemes, we follow papers~\cite{Muzychuk2012,Evdokimov2009a} and monograph~\cite{Zieschang2005a}. All the facts we use, can be found in these sources and references therein.

\subsection{Definitions.}
A pair $\cX=(\Omega,S)$, where $\Omega$ is a finite set and $S$ is a partition of $\Omega\times\Omega$, is called an {\it association scheme} or {\it scheme} on $\Omega$ if the following conditions are satisfied: $1_\Omega\in S$,  $S^*=S$, and given $r,s,t\in S$, the number
$$
c_{rs}^t:=|\alpha r\cap \beta s^*|
$$
does not depend on the choice of $(\alpha,\beta)\in t$.\medskip

The elements of $\Omega$, $S$, $S^\cup$, and the numbers~$c_{rs}^t$ are called the {\it points}, {\it basis relations}, {\it relations} and {\it intersection numbers} of~$\cX$, respectively. The numbers $|\Omega|$ and $|S|$ are called the {\it degree} and the {\it rank} of~$\cX$. A unique basic relation containing a pair $(\alpha,\beta)\in\Omega\times\Omega$ is denoted by $r(\alpha,\beta)$. Since the mapping $r:\Omega\times \Omega \to S$ depends only on $\cX$, it should be denoted by $r_\cX$, but we usually omit the subindex if this does not lead to confusion.

\subsection{Complex product.}
The set $S^\cup$ contains the relation $r\cdot s$ for all $r,s\in S^\cup$. It follows that this relation is the union (possibly empty) of basis relations of~$\cX$; the set of these relations is called the {\it complex product} of $r$ and $s$ and denoted by $rs$. Thus
$$
rs\subseteq S,\qquad r,s\in S.
$$
In what follows, for any $X,Y\subseteq S$, we denote by $XY$ the union of all sets $rs$ with $r\in X$ and $s\in Y$. Obviously, $(XY)Z=X(YZ)$ for all~$X,Y,Z\subseteq S$. 

\subsection{Valencies.}
For any basic relation $s\in S$, the number $|\alpha s|$ with $\alpha\in\Omega$ equals the intersection number $c_{ss^*}^{1}$, and hence does not depend on the choice of the point~$\alpha$. It is called the {\it valency} of $s$ and denoted by $n_s$; we say that $s$ is {\it thin} if $n_s=1$.\medskip

For the intersection numbers we have the following well-known identities:
\qtnl{181110a}
c_{r^*s^*}^{t^*}=c_{sr}^t\quad\text{and}\quad n_tc_{rs}^{t^*}=n_rc_{st}^{r^*}=n_sc_{tr}^{s^*},\qquad r,s,t\in S.
\eqtn

The scheme $\cX$ is said to be {\it regular} (respectively, {\it $\{1,k\}$-valenced}) if $n_s=1$ (respectively, $k>1$ and $n_s=1$ or $k$) for all $s\in S$. 

\subsection{Isomorphisms and schurity.}
A bijection from the point set of a scheme $\cX$ to the point set of a scheme~$\cX'$ is
called an {\it isomorphism} from $\cX$ to $\cX'$ if it induces a bijection between their sets of basis relations. The schemes $\cX$ and $\cX'$ are said to be {\it isomorphic} if there exists an isomorphism from $\cX$ to $\cX'$.\medskip

An isomorphism from a scheme $\cX$ to itself is called {\it automorphism} if the induced bijection on the basis relations of $\cX$ is the identity. The set of all automorphisms of a scheme~$\cX$ is a group with respect to composition and will be denoted by $\aut(\cX)$.\footnote{It is more natural to define an automorphism of a scheme as an isomorphism of it to itself; however, we should follow a long tradition in which a scheme is treated as a colored graph.}\medskip

Conversely, let $K\le\sym(\Omega)$ be a transitive permutation group, and let $S$ denote the set of orbits in the induced action of~$K$ on~$\Omega^2$. Then, $\cX=(\Omega,S)$ is a scheme; we say that $\cX$ {\it is  associated} with~$K$. A scheme on $\Omega$ is said to be {\it schurian} if it is associated with some transitive permutation group on~$\Omega$. A scheme~$\cX$ is schurian if and only if it is associated with the group~$\aut(\cX)$.

\subsection{Algebraic isomorphisms and separability.}
Let $\cX$ and $\cX'$ be schemes. A bijection $\varphi:S\to S',\ r\mapsto r'$ is called an {\it algebraic isomorphism} from~$\cX$ to~$\cX'$ if
\qtnl{f041103p1}
c_{rs}^t=c_{r's'}^{t'},\qquad r,s,t\in S.
\eqtn
In this case, $\cX$ and $\cX'$ are said to be {\it algebraically isomorphic}.\medskip

Each isomorphism~$f$ from~$\cX$ onto~$\cX'$ induces an algebraic isomorphism 
$$
\varphi_f:r\mapsto r^f
$$ 
between these schemes. The set of all isomorphisms inducing the algebraic isomorphism~$\varphi$ is denoted by $\iso(\cX,\cX',\varphi)$. In particular,
\qtnl{190316b}
\iso(\cX,\cX,\id_S)=\aut(\cX),
\eqtn
where $\id_S$ is the identity mapping on $S$. A scheme~$\cX$ is said to be {\it separable} if for any algebraic isomorphism~$\varphi:S\to S'$, the set $\iso(\cX,\cX',\varphi)$ is not empty.\medskip

The algebraic isomorphism $\varphi$ induces a bijection from $S^\cup$
onto $(S')^\cup$: the union $r\cup s\cup\cdots$ of basis relations
of $\cX$ is taken to $r'\cup s'\cup\cdots$. This bijection is also
denoted by $\varphi$. 

\subsection{Faithful maps.}
Let $\cX=(\Omega,S)$ and $\cX'=(\Omega',S')$ be schemes, and let $\varphi:S\to S'$ be
an algebraic isomorphism. A bijection $f$ from a subset of $\Omega$ to a subset of $\Omega'$ is said to be $\varphi$-\textit{faithful} if
$$
r(\alpha,\beta)^\varphi=r'(\alpha^f,\beta^f)\quad
\text{for all}\ \,\alpha,\beta\in\dom(f),
$$
where $\dom(f)$ is the domain of~$f$.  Clearly, if $f$ is a $\varphi$-faithful map, then the restriction of $f$ to any subset of~$\dom(f)$ is also $\varphi$-faithful.\medskip

A $\varphi$-faithful map $f$ is said to be {\it $\varphi$-extendable} to a point $\gamma\in\Omega$ if there exists a $\varphi$-faithful map with domain $\dom(f)\,\cup\,\{\gamma\}$, or, equivalently, if
\qtnl{100617u}
\bigcap_{\alpha\in\dom(f)}\alpha^fr(\alpha,\gamma)^\varphi\ne\varnothing.
\eqtn
A $\varphi$-faithful map which is $\varphi$-extendable  to every point of~$\Omega$, is said to be $\varphi$-extendable. From the definitions of schemes and algebraic isomorphisms, it follows that every $\varphi$-faithful map $f$ with $|\dom(f)|\le 2$ is $\varphi$-extendable. In these terms, one can give a sufficient condition for schurity and separability of a scheme.

\thrml{220718a}{\rm \cite[Corollary~2.2] {Hirasaka2018}}
Let $\cX$ be a scheme. Suppose that for every algebraic isomorphism $\varphi$ to another scheme, each $\varphi$-faithful map is $\varphi$-extendable. Then~$\cX$ is schurian and separable.
\ethrm

\section{Saturation condition}\label{281117l}

Throughout this section, $k>1$ is an integer  and  $\cX=(\Omega,S)$ is a scheme. By technical reasons, the basis relations of $\cX$ are mainly denoted below by $x,y,z$ rather than $r,s,t$. Our primary goal is to define a graph with vertex set
$$
S_k=\{x\in S:\ n_x=k\}
$$
that accumulates an information about intersection numbers equal to~$1$ (this graph was used implicitly in~\cite{Muzychuk2012a} an explicitly in~\cite{ChenPonom2017}). The following simple lemma (not formulated but proved in~\cite{ChenPonom2017}) indicates a way how we do this.

\lmml{230518p}
Given $x,y\in S_k$,
\qtnl{281117i}
|x^*y|=k\quad\Leftrightarrow\quad c_{xs}^y=1\ \,\text{for all}\ \,s\in x^*y.
\eqtn
\elmm
\proof We have $n_{x^{}}=n_{x^*}=n_{y^{}}=n_{y^*}=k$. By formulas~\eqref{181110a}, this implies that
$$
k^2=n_{x^*}n_{y^{}}=\sum_{s\in x^*y}n_sc_{x^*y^{}}^s=
\sum_{s\in x^*y}n_{y^*}c_{s^*x^*}^{y^*}=k\sum_{s\in x^*y}c_{xs}^y.
$$
Since $c_{xs}^y\ge 1$ for all $s\in x^*y$, we are done.\eprf\medskip

Let us define a relation $\sim$ on $S_k$ by setting $x\sim y$ if the right- or left-hand side in formula~\eqref{281117i}  holds true. This relation is symmetric, because 
$$
c_{xs}^y=\frac{n_{x^*}}{n_y}c_{sy^*}^{x^*}=c_{ys^*}^x
$$
for all $x,y\in S_k$. The (undirected) graph $\fX=\fX_k$ {\it associated} with the scheme~$\cX$ has vertex set $S_k$ and  adjacency relation~$\sim$. Note that this graph can contain loops.  

\dfntn
The scheme  $\cX$ is said to be {\it $k$-saturated} if for any set  $T\subseteq S$ with at most four elements, the set
\qtnl{280518d}
N(T)=\{y\in S_k:\ y\sim x\quad\text{for all}\ \,x\in T\}
\eqtn
is not empty. 
\edfntn

In this case, any two vertices of the graph $\fX$ are connected by a path of length at most two. A $k$-saturated $\{1,k\}$-scheme is said to be {\it saturated} and the mention of~$k$ is omitted. The following statement immediately follows from the definitions.

\lmml{250718a}
Any algebraic isomorphism from $\cX$ to $\cX'$ induces an isomorphism from the graph~$\fX$ to the graph~$\fX'$ associated with~$\cX'$. In particular, $\cX$ is $k$-saturated if and only if so is~$\cX'$.
\elmm

{\bf Example: schemes of affine spaces.} Let $\cA$ be a finite affine space with point set~$\Omega$ and line set~$L$, see~\cite{BuekenhoutCameron1995}. Denote by $\cP$ the set of parallel classes of lines. The lines belonging to a class $P\in\cP$  form a partition of $\Omega$; the corresponding equivalence relation with removed diagonal is denoted by~$e_P$. Using the axioms of affine spaces, one can easily verify that the set $S=S_\cA$ of all the~$e_P$ together with~$1_\Omega$ forms a commutative scheme such that
\qtnl{050816a}
c_{e_Pe_P}^s=\css
q-1       & $\text{if $s=1$,}$\\
q-2      & $\text{if $s=e_P$,}$\\
0          & $\text{otherwise,}$\\
\ecss
\eqtn
and
\qtnl{050816b}
c_{e_Pe_Q}^s=\css
1       & $\text{if $s\subseteq e_P\cdot e_Q$ and $s\not\in\{1,e_P,e_Q\}$,}$\\
0      & $\text{otherwise.}$\\
\ecss
\eqtn
where $q$ is the cardinality of a line. We say that $\cX=(\Omega,S)$ is the scheme associated with the affine space~$\cA$. Formulas~\eqref{050816a} and~\eqref{050816b} imply that $\cX$ is a  $\{1,q-1\}$-valenced scheme. Moreover, $\fX$ is a {\it complete} graph, i.e., any two distinct vertices form an edge.  In particular, the scheme~$\cX$ is $(q-1)$-saturated whenever the rank of~$\cX$ is at least~$6$.
\medskip

We complete the section by establishing  a sufficient condition for a scheme $\cX$ to be $k$-saturated in terms of  the indistinguishing number
$$
c(s)=\sum_{t\in S}c_{tt^*}^s
$$
of the relation~$s$, see~\cite{Muzychuk2012a}. One can see that this number is equal to 
the cardinality of the set
\qtnl{071217g}
\Omega_{\alpha,\beta}=\{\gamma\in\Omega:\ r(\gamma,\alpha)=r(\gamma,\beta)\}
\eqtn
for any $(\alpha,\beta)\in s$. 

\thrml{281117i1}\cite[Lemma~5.2]{ChenPonom2017}
Let $c$ be the maximum of $c(s)$, $s\in S^\#$. Then the scheme $\cX$  is $k$-saturated whenever 
$
|S_k|>4c(k-1).
$ 
\ethrm

\section{Desarguesian two-valenced schemes}\label{281117a4}

The concept of a Desarguesian scheme comes from a property of a geometry to be Desarguesian. Throughout this section, $k>1$ is an integer  and  $\cX=(\Omega,S)$ is a $\{1,k\}$-scheme. We also keep notation of Section~\ref{281117l}.\medskip

Let us define a noncommutative geometry associated with the scheme $\cX$ as follows: the points are elements of~$S$, the lines are the sets $x^*y$, $x,y\in S$, and the incidence relation is given by inclusion.  Thus the point $z\in S$ belongs to the line~$x^*y$ if and only if $z\in x^*y$. The geometry is extremely unusual: the line $x^*y$ does not necessarily contain the points $x,y$, and can be different from $y^*x$. However, in the terms of this geometry, one can define Desarguesian configurations, see below.\medskip

\begin{figure}[h]
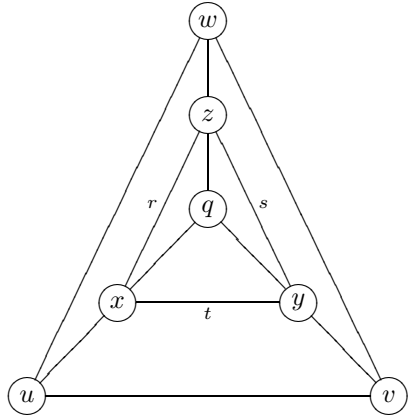

	\grphp{
		& & & & \VRT{w}\ar@{-}[dd] \ar@{-}[ddddddddllll] \ar@{-}[ddddddddrrrr]& &  & &\\
		& & & & & & & &\\
		& & & &\VRT{z}\ar@{-}[dd] \ar@{-}[ddddrr]^{s} \ar@{-}[ddddll]_{r}& &  & &\\
		& & & & & & & & \\
		& & &&\VRT{q}\ar@{-}[ddll] \ar@{-}[ddrr]& &  & &\\
		& & & &  & & & &\\
		& &\VRT{x} \ar@{-}[ddll] \ar@{-}[rrrr]_{t}&&&& \VRT{y} \ar@{-}[ddrr] & &\\
		& & & &  & & & &\\
		\VRT{u}  \ar@{-}[rrrrrrrr]& & & &  & & & & \VRT{v}\\
	}
	\caption{The Desarguesian configuration.}\label{281117a5}
\end{figure}

Assume that we are given two triangles with vertices $x,y,z\in S$ and $u,v,w\in S$, respectively, that are perspective with respect to a point~$q$, i.e.,
\qtnl{281117u1}
u\in x^*q,\quad v\in y^*q,\quad w\in z^*q,
\eqtn
see the configuration depicted in Fig.~\ref{281117a5}; note that the intersections of lines do not necessarily consist of a unique point, and even may be empty. However if
\qtnl{281117u7}
x^*z\cap uw^*=\{r\},\quad z^*y\cap wv^*=\{s\},\quad x^*y\cap uv^*=\{t\}
\eqtn
for some $r,s,t\in S$, then, as in the case of Desargues' theorem, we would like that the point $t$ would lie on the line $rs$. When this is true, this configuration is said to be {\it Desarguesian}. More precisely, the ten relations in Fig.~\ref{281117a5} form a Desarguesian configuration if conditions~\eqref{281117u1} and~\eqref{281117u7}  are satisfied and $t\in rs$.  In what follows we are going to study $\{1,k\}$-schemes with sufficiently many Desarguesian configurations.\medskip

Let $x,y,z\in S_k$ and $r,s\in S$  be basis relations of the scheme~$\cX$. We say that they form  an {\it initial} configuration if
\qtnl{240616t}
x\sim z\sim y\qaq r\in x^*z,\ \,s\in z^*y.
\eqtn
In geometric language, this means that the points $r$ and $s$ belong to the lines $x^*z$ and $z^*y$, respectively, and each of these lines consists of exactly~$k$ points. 

\begin{figure}[t]
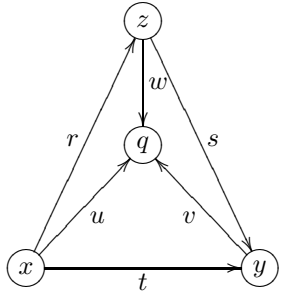

	\grphp{
		& & & \VRT{z}  \ar[ddd]^*{w}\ar[ddddddrrr]^*{s} & & & \\
		& & & & & & \\
		& & & & & & \\
		& & & \VRT{q} & & & \\
		& & & & & & \\
		& & & & & &\\
		\VRT{x} \ar[uuuuuurrr]^*{r}\ar[rrruuu]_*{u}\ar[rrrrrr]_*{t}
		& & &  & & &
		\VRT{y} \ar[llluuu]^*{v}\\
	}
	\caption{The relations $r$ and $s$ are linked with respect to $(x,y,z)$.}\label{f667}
\end{figure}

\dfntn
The relations $r$ and $s$ are said to be {\it linked} with respect to $(x,y,z)$  if  the initial configuration is contained in a Desarguesian configuration, namely, there exist
$$
q\in N(x,y,z),\qquad u,v,w\in S,\qquad t\in rs,
$$
for which conditions~\eqref{281117u1} and~\eqref{281117u7} are satisfied, where $N(x,y,z)=N(\{x,y,z\})$ (a more compact picture of the linked relations is given in Fig.\ref{f667}). 
\edfntn

Assume that the relations $r$ and $s$ are linked with respect to $(x,y,z)$.  Then the relation~$t$ is uniquely determined by the third of equalities~\eqref{281117u7}. The following statement shows that in this case, $t$ does not depend on the choice of~$q$ and~$u,v,w$. Below, we fix a point $\alpha\in\Omega$ and set 
\qtnl{190818o}
r_{x,y}=r\cap (\alpha x\times \alpha y)
\eqtn
for all $r\in S$ and $x,y\in S_k$. 

\lmml{090316u}
Assume that $r$ and $s$ are linked with respect to $(x,y,z)$.  Then
\qtnl{250518p}
r_{x,z}\cdot s_{z,y}\subseteq t_{x,y},
\eqtn
with equality if $x\sim y$.
\elmm
\proof By formulas~\eqref{281117u7}, we have
\qtnl{120316d}
u^{}_{x,q}\cdot w^*_{q,z}\subseteq r^{}_{x,z},\quad
w^{}_{z,q}\cdot v^*_{q,y}\subseteq s^{}_{z,y},\quad
u^{}_{x,q}\cdot  v^*_{q,y}\subseteq t^{}_{x,y}.
\eqtn
The relations $u^{}_{x,q}\cdot w^*_{q,z}$ and $r^{}_{x,z}$ are matchings, because $x\sim q\sim z$, and $x\sim z$. By the first inclusion in~\eqref{120316d}, this implies the first of the following two equalities, the second one is proved similarly: 
$$
u^{}_{x,q}\cdot w^*_{q,z}=r^{}_{x,z}\qaq w^{}_{z,q}\cdot v^*_{q,y}=s^{}_{z,y}.
$$ 
Now the third inclusion in~\eqref{120316d} yields
$$
r^{}_{x,z}\cdot s^{}_{z,y}=(u^{}_{x,q}\cdot w^*_{q,z})\cdot
(w^{}_{z,q}\cdot v^*_{q,y})=u^{}_{x,q}\cdot  v^*_{q,y}
\subseteq t^{}_{x,y},
$$
which proves~\eqref{250518p}. If $x\sim y$, then the relation $t_{x,y}$ is a matching and we are done by the above argument.\eprf\medskip

At this  point, we need an auxiliary statement. In the geometric language, the conclusion of this statement means that the lines~$rs$ and~$x^*y$ have a unique common point.

\lmml{270518a1}{\rm \cite[Theorem~5.1]{Muzychuk2012a}} 
Let $x,y,z\in S_k$ be such that
\qtnl{291117a}
(xx^*\,yy^*)\cap zz^*=\{1\}.
\eqtn 
Then  $|rs\cap x^*y|=1$ for all $r\in x^*z$ and $s\in z^*y $.
\elmm

The statement below establishes two  sufficient conditions for  relations $r$ and $s$ to be linked with respect to $(x,y,z)$. 

\crllrl{271117u1}
Let $x,y,z\in S_k$ and $r,s\in S$ form an initial configuration. Then~$r$ and~$s$ are linked with respect to $(x,y,z)$, whenever at least one of the following statements holds:
\nmrt
\tm{L1} $z$ is a loop of the graph~$\fX$ and condition~\eqref{291117a} is satisfied,
\tm{L2} there exists $q\in S_k$ such that
\qtnl{020816c}
qq^*\,\cap\,(xx^*yy^*\,\cup\,xx^*zz^*\,\cup\,zz^*yy^*)=\{1\}.
\eqtn
\enmrt
\ecrllr
\proof Assume that (L1) holds. Then by Lemma~\ref{270518a1}, we have $rs\cap x^*y=\{t\}$ for some $t\in S$. Then condition \eqref{281117u7} is obviously satisfied for $q=z$ and $u=r$, $v=s^*$, $w=1$.\medskip 

Now assume that (L2) holds. Fix a point $\alpha$. Since $r\in x^* z$ and $s\in z^* y$, one can find points $\beta\in \alpha x$, $\gamma\in \alpha z$, and $\delta\in \alpha y$ such that 
\qtnl{201118a}
(\beta,\gamma)\in r\qaq (\gamma,\delta)\in s.
\eqtn
Now  take $\epsilon \in \alpha q$ and set
$$
u:=r(\beta,\epsilon),\quad v:=r(\delta,\epsilon),\quad w:=r(\gamma,\epsilon).
$$
Then formula~\eqref{281117u1} holds. Moreover by the condition~(L2), the hypothesis of Lemma~\ref{270518a1} is satisfied in the following three cases:
\begin{align*}
(x,y,z)&=(x,z,q),\  (r,s)=(u,w^*),\\
(x,y,z)&=(z,y,q),\  (r,s)=(w,v^*), \\
(x,y,z)&=(x,y,q),\  (r,s)=(u,v^*).
\end{align*}
Therefore formula ~\eqref{281117u7} also holds with $t=r(\beta,\gamma)$.  Finally, $t\in rs$ in view of~\eqref{201118a}.  It follows that the configuration formed by $x,y,z$, $u,v,w$, $r,s,t$, and $q$  is Desarguesin. Thus $r$ and~$s$ are linked with respect to $(x,y,z)$.\eprf\medskip

Now we arrive to the main definition in this section. Namely, the scheme~$\cX$ is said to be {\it Desarguesian} with respect to $S_k$ if  for all $x,y,z\in S_k$ and all $r,s\in S$  satisfying \eqref{240616t}, the elements $r$ and $s$ are linked with respect to $(x,y,z)$.  When the scheme is two-valenced, the mention of $S_k$ is omitted. The following statement immediately follows from the definitions.

\lmml{250718a1}
Let $\cX$ and $\cX'$ be algebraically isomorphic two-valenced schemes. Then~$\cX$ is Desarguesian if and only if so is~$\cX'$.
\elmm

{\bf Example: schemes of affine spaces (continuation).}
Let $\cX$ be the scheme associated with an affine space $\cA$ of order~$q$ and dimension  at least~$3$. Then from formulas~\eqref{050816a} and~\eqref{050816b} it follows that for any three  parallel classes $P$,  $Q$, and $R$ there exists a parallel class~$T$ such that
$$
e_T\not\in(e_Pe_Q\,\cup\,e_Pe_R\,\cup\,e_Re_Q).
$$
It follows that the statement (L2) of Corollary~\ref{271117u1} holds for $q=e_T$, $x=e_P$, $y=e_Q$, and $z=e_R$. Therefore any relations $r\in e_Pe_R$ and $s\in e_Re_Q$ are linked with respect to $(e_P,e_Q,e_R)$. Thus, the scheme $\cX$ is Desarguesian. If the space~$\cA$ is an affine plane, i.e., an affine space of dimension~$2$, then 
$\cX$ is not Desarguesian.\medskip 

\section{Proof of Theorem~\ref{181117a}}\label{090818n}

Let $\cX=(\Omega,S)$, $\cX'=(\Omega',S')$, and $\varphi$ an algebraic isomorphism from~$\cX$ onto~$\cX'$. Then obviously $\cX'$ is two-valenced. Moreover, $\cX'$ is saturated and Desarguesian by Lemmas~\ref{250718a} and~\ref{250718a1}, respectively. By Theorem~\ref{220718a}, it suffices to verify that given $\alpha,\beta\in\Omega$,  any $\varphi$-faithful map 
$$
f:\{\alpha,\beta\}\to \Omega'
$$ 
is extended to a combinatorial isomorphism $\wh f:\cX\to\cX'$. Without loss of generality, we may assume that $r(\alpha,\beta)\in S_k$, where $k>1$ is the valency of~$\cX$ (and $\cX'$). Indeed, in this case the isomorphism $\wh f$ extends any faithful map $f:\{\alpha,\wt\beta\}\to \Omega'$, where $r(\alpha,\wt\beta)\in S_1$. \medskip

By the definition of the graph $\fX=\fX_k$, the union of the sets 
\begin{align*}
\Delta_0&:=\alpha S_1\\
\Delta_1& :=\{\delta\in\alpha S_k:\ r(\alpha,\beta)\sim r(\alpha,\delta)\}\\
\Delta_2&:=\{ \delta\in\alpha S_k:\ r(\alpha,\beta)\not\sim r(\alpha,\delta)\}
\end{align*}
equals $\Omega$. Clearly, $\beta\in\Delta_1$ if and only if  $r(\alpha,\beta)$ forms a loop of the graph $\fX$. The statement below immediately follows from the definition of the relation~$\sim$, Lemma~\ref{250718a}, and the fact that $\fX$ is of diameter at most~$2$ (the saturation condition).

\lmml{250718b}
For any point $\delta\in\Omega$, the following statements hold:
\nmrt
\tm{1} if $\delta\in\Delta_0$, then the relation $r(\alpha,\delta)$ is thin,
\tm{2}  if $\delta\in\Delta_1$, then $c_{r(\alpha,\delta),r(\delta,\beta)}^{r(\alpha,\beta)}=1$,
\tm{3} if $\delta\in\Delta_2$, then there exists $\gamma\in\Delta_1$ such that
$c_{r(\alpha,\gamma),r(\gamma,\delta)}^{r(\alpha,\delta)}=1$.
\enmrt
\elmm

By Lemma~\ref{250718b}, one can define a mapping $\wh f:\Omega\to \Omega'$ by the following conditions:
$$
\{\delta^{\wh f}\}=\css
\alpha^f r(\alpha,\delta)^\varphi   &\text{ if $\delta\in \Delta_0$},\\
\alpha^f r(\alpha,\delta)^\varphi \cap \beta^f  r(\beta,\delta)^\varphi &\text{ if $\delta\in \Delta_1$},\\
\alpha^f r(\alpha,\delta)^\varphi \cap \gamma^{\wh f} r(\gamma,\delta)^\varphi
&\text{ if $\delta\in \Delta_2$},\\
\ecss
$$
where in the last line, $\gamma$ is as in statement~(3) of Lemma~\ref{250718b}, i.e.,
$$
r(\alpha,\beta) \sim r(\alpha,\gamma) \sim r(\alpha,\delta).
$$
From the definition of the graph~$\fX$, it follows that the $\wh f$-image is uniquely determined for all $\delta\in \Delta_0\cup \Delta_1$. However, if $\delta\in\Delta_2$, then the image depends on the choice of the point~$\gamma$.\medskip

We complete the proof by verifying  that  the mapping  ${\wh f}$ is  $\varphi$-faithful, i.e.,  	
\qtnl{230718a}
r(\delta,\epsilon)^\varphi=r(\delta^{\wh f},\epsilon^{\wh f})
\eqtn
for all $\delta,\epsilon \in \Omega$. Note that this is true by the definition of $\wh f$ whenever $\alpha\in\{\delta,\epsilon\}$. In what follows we need the following lemma.

\lmml{051118a}
Let $\delta,\gamma,\epsilon\in\Omega$ be such that
\nmrt
\tm{1} $r(\alpha,\delta)\sim r(\alpha,\gamma)\sim r(\alpha,\epsilon)$,
\tm{2} $\wh f|_{\{\delta,\gamma\}}$ and $\wh f|_{\{\gamma,\epsilon\}}$ are $\varphi$-faithful.
\enmrt
Then $\wh f|_{\{\delta,\epsilon\}}$ is $\varphi$-faithful.
\elmm
\proof Let $x=r(\alpha,\delta)$, $y=r(\alpha,\epsilon)$, and $z=r(\alpha,\gamma)$. Then by condition~(1), the triple $(x,y,z)$ together with relations $r=r(\delta,\gamma)$, $s=r(\gamma,\epsilon)$ forms an initial configuration. The Desarguesian condition implies that $r$ and $s$ are linked with respect to $(x,y,z)$.  Therefore the relation $t=r(\delta,\epsilon)$ is uniquely determined by the relations $r_{x,z}$ and~$r_{z,y}$ (Lemma~\ref{090316u}). By condition~(2), we have
$$
(r_{x,z})^{\wh f}=(r^\varphi)_{x^\varphi, z^\varphi}\qaq
(s_{z,y})^{\wh f}=(s^\varphi)_{z^\varphi, y^\varphi}
$$
and hence
$$
r(\delta^{\wh f},\epsilon^{\wh f})\in r(\delta^{\wh f},\gamma^{\wh f})r(\gamma^{\wh f},\epsilon^{\wh f})=r(\delta,\gamma)^\varphi r(\gamma,\epsilon)^\varphi.
$$ 
Since also
$$
r(\delta,\epsilon)^\varphi\in (r(\delta,\gamma)r(\gamma,\epsilon))^\varphi=r(\delta,\gamma)^\varphi r(\gamma,\epsilon)^\varphi,
$$
we conclude that  
$$
r(\delta^{\wh f},\epsilon^{\wh f})=r(\delta,\epsilon)^\varphi,
$$
as required.\eprf\medskip

The rest of the proof of formula~\eqref{230718a} is divided for several cases considered separately.\medskip
 	
{\bf Case 1:}  $\delta \in \Delta_0$ and $\epsilon \in \Omega$. Here, the relations $r(\delta,\alpha)$ and $r(\delta^{\wh f},\alpha^{\wh f})$ are thin. Therefore, 
$$
\{r(\delta,\epsilon)\}=r(\delta,\alpha)r(\alpha,\epsilon)\qaq 
r(\delta^{\wh f},\alpha^{\wh f})r(\alpha^{\wh f},\epsilon^{\wh f})=\{r(\delta^{\wh f},\epsilon^{\wh f})\}.
$$
This implies that
$$
\{r(\delta,\epsilon)^\varphi\}=r(\delta,\alpha)^\varphi r(\alpha,\epsilon)^\varphi=r(\delta^{\wh f},\alpha^{\wh f})r(\alpha^{\wh f},\epsilon^{\wh f})=\{r(\delta^{\wh f},\epsilon^{\wh f})\},
$$ 	
as required.\medskip 

{\bf Case 2:} $\delta,\epsilon \in \Delta_1$. 	The required statement follows from Lemma~\ref{051118a} for $\gamma=\beta$.\medskip

\begin{figure}[t]
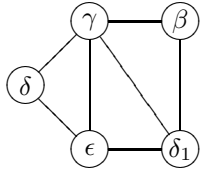

\grphp{
 & \VRT{\gamma}  \ar@{-}[rr]   \ar@{-}[dd]  \ar@{-}[ddrr]& &\VRT{\beta} \ar@{-}[dd]\\
\VRT{\delta}\ar@{-}[ru]	\ar@{-}[rd]	& &  &\\
 & \VRT{\epsilon}  \ar@{-}[rr] & &\VRT{\delta_1}\\
}
	\caption{Configuration for Case~3.}\label{260718a}
\end{figure}

{\bf Case 3:}  $\delta \in \Delta_1$ and $\epsilon \in \Delta_2$.  	By the definition of $\Delta_2$ there exists $\delta_1\in \Delta_1$ such that
$$
r_\beta\sim r_{\delta_1}\sim r_\epsilon\qaq \{\epsilon^{\wh f}\}=\alpha^f (r_\epsilon)^\varphi \cap \delta_1^{\wh f} r(\delta_1,\epsilon)^\varphi,
$$
where $r_\beta=r(\alpha,\beta)$, $r_{\delta_1}=r(\alpha,\delta_1)$, and $r_\epsilon=r(\alpha,\epsilon)$. 	By the saturated condition, there exists $\gamma\in \Delta_1$ such that 
$$
r_\gamma\in N(r_\beta,r_\delta,r_{\delta_1},r_\epsilon),
$$
the obtained configuration is depicted in Fig.~\ref{260718a}.
From Case~$2$, it follows that the conditions of Lemma~\ref{051118a} are satisfied for $$
(\delta,\gamma,\epsilon)=(\gamma,\delta_1,\epsilon).
$$ 
Therefore $\wh f|_{\{\epsilon,\gamma\}}$ is $\varphi$-faithful. This shows that  the conditions of Lemma~\ref{051118a} are satisfied and hence  $\wh f|_{\{\delta,\epsilon\}}$ is $\varphi$-faithful.\medskip

{\bf Case 4:}  $\delta, \epsilon \in \Delta_2$.  By the definition of $\Delta_2$ there exists $\delta_1\in \Delta_1$ such that
$$
r_\beta\sim r_{\delta_1}\sim r_\delta\qaq \{\delta^{\wh f}\}=\alpha^f (r_\delta)^\varphi \cap \delta_1^{\wh f} r(\delta_1,\delta)^\varphi.
$$
By the saturation condition, there exists $\gamma\in\Delta_1$ such that
$$
r_\gamma\in N(r_\beta,r_\delta, r_{\delta_1}, r_\epsilon).
$$
By the similar argument as in Case~3, we obtain $r(\delta,\epsilon)^\varphi=r(\delta^{\wh f},\epsilon^{\wh f})$, where Case~3 is used to deal with $r(\delta,\delta_1)$, $r(\delta,\gamma)$, and $r(\epsilon,\gamma)$.\eprf\medskip

\section{Proofs of corollaries and Theorem~\ref{thm111}}\label{090818a}
	
{\bf Proof of Corollary~\ref{110317a}}	
Set $f(x)=3x^6+1$. By Theorem~\ref{181117a}, it suffices to verify that any pseudocyclic scheme $\cX$  of degree $n$ and valency~$k$,  is saturated and Desarguesian, whenever
\qtnl{071217e}
n>3k^6.
\eqtn

By \cite[Theorem~3.2]{Muzychuk2012a}, we have  $n_s=k$ and $c(s)=k-1$ for any irreflexive basis relation $s$ of~$\cX$.  It follows that for $k\ge 2$, 
$$
|S_k|=\frac{n-1}{k}> \frac{3k^6-1}{k}> 4k(k-1)=4ck.
$$
Thus, $\cX$ is saturated by Theorem~\ref{281117i1}.\medskip

To prove that $\cX$ is Desarguesian, let $x,y,z\in S_k$ (here $S_k=S^\#$ and we do not assume that $x\sim z\sim y$).  By Corollary~\ref{271117u1}, it suffices to find $q\in S_k$ such that condition~\eqref{020816c} is satisfied. Assume on the contrary that no $q$ satisfies this condition. Then for a fixed $\alpha\in\Omega$ and each $q\in S_k$, there exists 
$$
\beta_q\in\alpha \underbrace{(xx^*yy^*\,\cup\,xx^*zz^*\,\cup\,zz^*yy^*)}_{T}
$$ 
other than $\alpha$ and such that $r(\alpha,\beta_q)\in qq^*$.  It follows that
\qtnl{071217i}
\Omega_{\alpha,\beta_q}\,\cap\, \alpha q^*\ne\varnothing,
\eqtn
where the set $\Omega_{\alpha,\beta_q}$ is as in \eqref{071217g}. Since 
$$
|\alpha T|\le n_x^2n_y^2+n_y^2n_z^2+n_z^2n_x^2=3k^4
\qaq
|S_k|=\frac{n-1}{k}\ge \frac{3k^6}{k}=3k^5,
$$
there exist a point $\beta\ne\alpha$ such that $\beta=\beta_q$ for at least $k$ relations~$q\in S_2$. In view of \eqref{071217i}, this implies that for $s=r(\alpha,\beta)$ we have
$$
k-1=c(s)\ge |\Omega_{\alpha,\beta}|\ge k,
$$
a contradiction.\eprf\medskip

{\bf Proof of Corollary~\ref{280618a}.} Let $\cX=(\Omega,S)$ be a quasi-thin scheme satisfying the condition $n_{s^*s}\ne 2$ for all $s\in S$. Without loss of generality, we may assume that $|\Omega|> 24$.\medskip

{\bf Claim.} The graph $\fX$ is complete and has at least~$9$ vertices. Moreover, for any  $x,y\in S_2$,
\qtnl{030818a}
xx^*=yy^*\quad\Leftrightarrow\quad x=y.
\eqtn

\proof  From \cite[Lemma~4.1]{Hirasaka2001}, it follows that given $x,y\in S_2$ and $s\in S$, the number $c_{xs}^y$ is at most~$2$, and by \cite[Lemma~3.1]{Hirasaka2002a} also
$$
c_{xs}^y=2\quad\Leftrightarrow\quad s\in S_2,\ ss^*\subset S_1,
$$
Since $n_{ss^*}\ne 2$, this implies \eqref{030818a} and shows that in our case, $c_{xs}^y\le 1$ for all $x,y\in S_2$; in particular, $\fX$ is a complete graph with loops. Furthermore,  again by   \cite[Lemma~3.1]{Hirasaka2002a}, in our case $xs\not\subset S_1$ for all $x,s\in S_2$. It follows that 
$$
|\Omega|\ge 3|S_1|.
$$ 
Assume on the contrary that $\fX$ has at most~$8$ vertices, i.e., $|S_2|\le 8$. Then 
$$
3|S_1|\le|\Omega|=|S_1|+2|S_2|\le |S_1|+16,
$$
whence $|S_1|\le 8$. It follows that $|\Omega|\le 24$, a contradiction.\eprf\medskip

 By Theorem~\ref{181117a}, it suffices to verify that $\cX$ is saturated and Desarguesian. The former follows from the Claim. To prove that $\cX$ is Desarguesian let $x,y,z\in S_k$ and $r,s\in S$  be basis relations of~$\cX$  that form an initial configuration. We have to verify that $r$ and $s$ are linked with respect to $(x,y,z)$.  Note that since $n_z=2$, the vertex $z$ forms  a loop of the graph~$\fX$. Therefore if  $(xx^*\,yy^*)\cap zz^*=\{1\}$, then we are done by Corollary~\ref{271117u1} (see the statement~(L1)).\medskip
 
 Let $zz^*\subseteq xx^*\cap yy^*$. Then the intersection of $S_2$ with the set
$$
xx^*yy^*\,\cup\,xx^*zz^*\,\cup\,zz^*yy^*=\{1,x^\ddp,y^\ddp\}\,\cup\,x^\ddp y^\ddp\,\cup\,x^\ddp z^\ddp\,\cup\,z^\ddp y^\ddp
$$
contains at most $8$ elements,  where for any $s\in S_2$, we set $s^\ddp$ to be  the unique non-thin element of~$ss^*$. On the other hand, in view of the Claim, the graph $\fX$ contains at least $9$ vertices. Thus there exists  $q'\in S_2$ such that
$$
q'\not\in xx^*yy^*\,\cup\,xx^*zz^*\,\cup\,zz^*yy^*.
$$
In view of formula~\eqref {030818a} and the assumption $n_{ss^*}\ne 2$, one can find $q\in S_2$, for which $q'=q^\ddp$. Now the required statement follows from Corollary~\ref{271117u1} (see the statement~(L2)).\eprf\medskip

{\bf Proof of Theorem~\ref{thm111}.} Let $\cX=(\Omega,S)$ be a meta-thin $\{1,p\}$-scheme. Then the group formed by the thin basis relations of $\cX$ (with respect to the composition) contains a subgroup (the thin residue) generated by the sets $ss^*$, $s\in S$.   Assume that this subgroup is isomorphic to an elementary abelian $p$-group of rank greater  than two. By Theorem~\ref{181117a}, it suffices to verify that $\cX$ is saturated and Desarguesian. \medskip

To prove that $\cX$ is saturated, let $s_i\in S_p$, $1\le i\le 4$. By the assumption there exist $u_1,u_2, u_3\in S$ such that
$$
\grp{u^{}_1u_1^*,u^{}_2u_2^*,u^{}_3u^*_3}\simeq E_{p^3}.
$$
Denote by $T$ the set $\{s^{}_1s^*_1,\ldots, s^{}_4s^*_4\}$. First, assume that $u^{}_ju_j^*\notin  T$ for some $j$. Then
$$
u^{}_ju_j^*\cap s^{}_is_i^\ast=\{1\},\qquad 1\le i\le 4.
$$
In view of \cite[Lemma~2.3]{Muzychuk2012a}, this implies that $u_j$ is adjacent in the graph $\fX$ with each of the $s_i$.\medskip

Now without loss of generality we may assume that  $u^{}_ju_j^*\in  T$ for all~$j$. Since $u^{}_ju_j^*\cap u^{}_ku_k^*=\{1\}$ for all distinct $j,k$, there exist  $i\in \{1,2,3,4\}$ and $j\in\{1,2,3\}$ such that 
$$
s^{}_is_i^*=u^{}_ju_j^*\ne s^{}_ks_k^*,\quad k\ne i.
$$
Taking into account that every element in $S_p$ has a loop, we conclude that $s_i$ is adjacent in the graph $\fX$ with each of the $s_i$ with $i\in \{1,2,3,4\}$. Thus the scheme~$\cX$ is saturated.\medskip

To prove that the scheme $\cX$ is Desarguesian, let $x,y,z\in S_k$ and $r,s\in S$  be basis relations of~$\cX$  that form an initial configuration. We have to verify that $r$ and~$s$ are linked with respect to $(x,y,z)$.  If 
$$
\grp{xx^*, yy^*, zz^*}\simeq E_{p^3},
$$
then this follows from Corollary~\ref{271117u1} (see the statement~(L1)). Assume that this condition is not satisfied, i.e., the group on the left-side is isomorphic to $E_p$ or $E_{p^2}$.  Then by the assumption of the theorem, there exists $q\in S_p$ such that
$$
qq^*\cap \grp{xx^*, yy^*,zz^*}=\{1\}.
$$
It follows that condition~\eqref{020816c} is satisfied. Thus the required statement follows from Corollary~\ref{271117u1} (see the statement~(L2)).\eprf

\end{document}